# ON A CLASS OF SELF-SIMILAR 2D SURFACE WATER WAVES

SIJUE WU

ABSTRACT. We construct a class of self-similar surface water waves and study its properties. This class of surface waves appears to be in very good agreement with a common type of wave crests in the ocean.

## 1. INTRODUCTION

The focus of this paper is on understanding singularities of the self-similar type in 2-D surface water waves.

The mathematical problem of $n$-dimensional water wave concerns the motion of the interface separating an inviscid, incompressible, irrotational fluid, under the influence of gravity, from a region of zero density (i.e. air) in $n$-dimensional space. It is assumed that the fluid region is below the air region. Assume that the density of the fluid is 1, the gravitational field is $-\mathfrak{g}\mathbf{k}$, where $\mathfrak{g} > 0$, $\mathbf{k}$ is the unit vector pointing in the upward vertical direction, and at time $t \geq 0$, the free interface is $\Sigma(t)$, and the fluid occupies region $\Omega(t)$. When surface tension is zero, the motion of the fluid is described by

$$\begin{cases} \mathbf{v}_t + \mathbf{v} \cdot \nabla \mathbf{v} = -\mathfrak{g}\mathbf{k} - \nabla P & \text{on } \Omega(t),\ t \geq 0, \\ \operatorname{div} \mathbf{v} = 0, \quad \operatorname{curl} \mathbf{v} = 0, & \text{on } \Omega(t),\ t \geq 0, \\ P = 0, & \text{on } \Sigma(t), \end{cases} \quad (1.1)$$

where $\mathbf{v}$ is the fluid velocity, $P$ is the fluid pressure. It is well-known that when surface tension is neglected, the water wave motion can be subject to the Taylor instability [6, 27, 5]. Assume that the free interface $\Sigma(t)$ is described by $z = z(\alpha, t)$, where $\alpha \in R^{n-1}$ is the Lagrangian coordinate, i.e. $z_t(\alpha, t) = \mathbf{v}(z(\alpha, t), t)$ is the fluid velocity on the interface, $z_{tt}(\alpha, t) = (\mathbf{v}_t + \mathbf{v} \cdot \nabla \mathbf{v})(z(\alpha, t), t)$ is the acceleration. Let $\mathbf{n}$ be the unit normal pointing out of $\Omega(t)$. The Taylor sign condition relating to Taylor instability is

$$-\frac{\partial P}{\partial \mathbf{n}} = (z_{tt} + \mathfrak{g}\mathbf{k}) \cdot \mathbf{n} \geq c_0 > 0 \quad (1.2)$$

point-wisely on the interface for some positive constant $c_0$. In [29, 30], we showed that the Taylor sign condition (1.2) always holds for the $n$-dimensional infinite depth water wave problem (1.1), $n \geq 2$, as long as the interface is non-self-intersecting; and the initial value problem of the water wave system (1.1) is uniquely solvable locally in time in Sobolev

Financial support in part by NSF grants DMS-0800194, DMS-1101434.

spaces for arbitrary given data. Earlier work includes Nalimov [22] and Yosihara [33] on local existence and uniqueness for small data in 2D. In [3, 12, 24, 28], rigorous justifications of the KdV, KP, Boussinesq, shallow water, NLS and various other asymptotic models from the full water wave equations were obtained, establishing in rigorous mathematical terms the typical surface wave behaviors in corresponding regimes. In [31, 32], we proved that the nature of the nonlinearity of the full water wave equation (1.1) is of cubic or higher orders, we then showed that for data that are smooth and small in some generalized Sobolev spaces, solutions of equation (1.1) remain small and smooth for almost global time in 2-D, and for all time in 3-D. In [16], for the 3-D water wave equation (1.1), the authors constructed a different class of smooth and small data, and showed that for these data the solutions remain small and smooth for all time and scatter. There are much more work on water waves recently [2, 9, 10, 17, 20, 21, 23, 25, 34]. Among them we mention the existence of the so called splash singularities from smooth data for the 2-D and 3-D water waves [8, 7, 11]; and the local wellposedness in some low regularity Sobolev classes for the Cauchy problem of the gravity waves [1].

Self-similarity is an important tool in mathematical physics and in the study of singularities [4, 14]. In situations in which no explicit reference length appear, such as when the boundary effect and external forces are negligible, typical phenomenon are often self-similar. Self-similarity has been exploited for the water waves. In [26], a class of self-similar solutions for a linear water wave equation was constructed; and in [15], it was done for a nonlinear approximation of the water wave equation (1.1).

Indeed, for the full water wave equation (1.1), a similarity law holds. An even wider set of similarity laws hold if the gravity is neglected (i.e. $\mathfrak{g} = 0$). In the next section, we will discuss the similarity laws for water waves with or without surface tension or gravity. In subsequent sections, we will construct a class of self-similar surface waves and study their properties. The self-similar solutions we construct appear to be in very good agreement with a common type of wave phenomena we observe (see 3).

In what follows we will focus on the two dimensional water waves.

## 2. Similarity Laws

We use complex variables and identify $z = (x, y)$ with $z = x + i y$. $\bar{z}$ is the complex conjugate of $z$; $\text{Im } z$, $\text{Re } z$ are respectively the imaginary and real parts of $z$. In this section, we discuss similarity laws for water waves with or without surface tension or gravity.

Let $z = z(\alpha, t)$, $\alpha \in R$ be the interface $\Sigma(t)$ at time $t$ in Lagrangian coordinate $\alpha$. We know its curvature is given by $\text{Im} \frac{\bar{z}_\alpha z_{\alpha\alpha}}{|z_\alpha|^3}$. Let $\sigma$ be the surface tension coefficient. The equation describing water waves with surface tension is given by

$$\begin{cases} \mathbf{v}_t + \mathbf{v} \cdot \nabla \mathbf{v} = -\mathfrak{g}\,i - \nabla P & \text{on } \Omega(t),\ t \geq 0, \\ \text{div}\,\mathbf{v} = 0, \quad \text{curl}\,\mathbf{v} = 0, & \text{on } \Omega(t),\ t \geq 0, \\ P = \sigma \text{Im} \frac{\bar{z}_\alpha z_{\alpha\alpha}}{|z_\alpha|^3}, & \text{on } \Sigma(t), \end{cases} \quad (2.1)$$

where $\mathbf{v}$ is the fluid velocity, $P$ is the fluid pressure. Notice that on the interface,

$$\bar{z}_\alpha \nabla P(z(\alpha,t),t) = \partial_\alpha\{P(z(\alpha,t),t)\} + i|z_\alpha|\frac{\partial P}{\partial \mathbf{n}}.$$

Similar to [29, 31], we can rewrite equation (2.1) into the following equivalent nonlinear evolution equation of the interface $z = z(\alpha, t)$, $\alpha \in R$:

$$\begin{cases} z_{tt} + \mathfrak{g}\,i = -\sigma \frac{z_\alpha}{|z_\alpha|^2}\partial_\alpha(\text{Im}\frac{\bar{z}_\alpha z_{\alpha\alpha}}{|z_\alpha|^3}) + i\,\mathfrak{a} z_\alpha \\ \bar{z}_t \text{ is the boundary value of a holomorphic function on } \Omega(t) \end{cases} \quad (2.2)$$

where $\mathfrak{a} = -\frac{1}{|z_\alpha|}\frac{\partial P}{\partial \mathbf{n}}$. Notice that $\bar{z}_t(\alpha, t) = \bar{\mathbf{v}}(z(\alpha,t),t)$. It is easy to check the following similarity laws for (2.1) and (2.2):

If $\mathbf{v} = \mathbf{v}(z,t)$, $P = P(z,t)$ with interface $\Sigma(t) : z = z(\alpha,t)$ are solutions for equations (2.1) and (2.2),

Case 1: if $\mathfrak{g} = \sigma = 0$, then for all $s$, all $\lambda > 0$,

$$\mathbf{v}_\lambda = \lambda^{s-1}\mathbf{v}(\lambda z, \lambda^s t), \quad P_\lambda = \lambda^{2(s-1)} P(\lambda z, \lambda^s t) \quad (2.3)$$

$$z_\lambda = \lambda^{-1} z(\lambda\alpha, \lambda^s t), \quad \mathfrak{a}_\lambda = \lambda^{2s-1}\mathfrak{a}(\lambda\alpha, \lambda^s t) \quad (2.4)$$

are also solutions of (2.1) and (2.2),

Case 2: if $\mathfrak{g} \neq 0$, $\sigma = 0$, then the similarity laws (2.3) (2.4) hold for $s = 1/2$.

Case 3: if $\mathfrak{g} = 0$, $\sigma \neq 0$, then the similarity laws (2.3) (2.4) hold for $s = 3/2$.

In the case when both $\mathfrak{g} \neq 0$, $\sigma \neq 0$, there are no similarity laws.

Now let $z = z(\alpha, t)$ be a solution of (2.2). Assume that a singularity occurs at $(\alpha, t) = (0, 0)$ in a self-similar way. We blow-up $z(\alpha, t)$ by letting $z^\epsilon(\alpha, t) = \epsilon^{-1} z(\epsilon\alpha, \epsilon^s t)$, and assume that as $\epsilon \to 0$ the limiting profile $Z = Z(\alpha, t)$ exists:

$$z^\epsilon(\alpha, t) = \epsilon^{-1} z(\epsilon\alpha, \epsilon^s t) \to Z(\alpha, t), \quad \text{as } \epsilon \to 0$$

and derivatives of $z^\epsilon$ approach the corresponding derivatives of $Z$. It is clear that $Z = Z(\alpha, t)$ is self-similar:

$$Z(\alpha, t) = \lambda^{-1} Z(\lambda\alpha, \lambda^s t), \quad Z_t(\alpha, t) = \lambda^{s-1} Z_t(\lambda\alpha, \lambda^s t) \quad \text{for all } \lambda > 0 \quad (2.5)$$

Let $\mathfrak{a}^\epsilon(\alpha, t) = \epsilon^{2s-1}\mathfrak{a}(\epsilon\alpha, \epsilon^s t)$ and assume

$$\mathfrak{a}^\epsilon(\alpha, t) \to A(\alpha, t) \qquad \text{as } \epsilon \to 0.$$

From the assumption that $z = z(\alpha, t)$ satisfies (2.2), we have

$$\begin{cases} \partial_t^2 z^\epsilon + \epsilon^{2s-1}\mathfrak{g}\, i = i\mathfrak{a}^\epsilon \partial_\alpha z^\epsilon - \sigma\epsilon^{2s-3}\dfrac{\partial_\alpha z^\epsilon}{|\partial_\alpha z^\epsilon|^2}\partial_\alpha(\mathrm{Im}\,\dfrac{\partial_\alpha \bar{z}^\epsilon \partial_\alpha^2 z^\epsilon}{|\partial_\alpha z^\epsilon|^3}) \\ \bar{z}_t^\epsilon \text{ is the boundary value of a holomorphic function on } \Omega^\epsilon(t) \end{cases} \qquad (2.6)$$

where $\Omega^\epsilon(t)$ is the rescaled fluid domain with boundary $\Sigma^\epsilon(t): z^\epsilon = z^\epsilon(\alpha, t)$. We see that if $\sigma = 0$, the gravity and the terms $\partial_t^2 z^\epsilon$ and $i\mathfrak{a}^\epsilon \partial_\alpha z^\epsilon$ achieve a balance when $s = 1/2$. If $\mathfrak{g} = 0$, the surface tension term achieves a balance with terms $\partial_t^2 z^\epsilon$ and $i\mathfrak{a}^\epsilon \partial_\alpha z^\epsilon$ when $s = 3/2$. If both surface tension and gravity are neglected, then $Z(\alpha, t)$ satisfies

$$\begin{cases} Z_{tt} = i\, A Z_\alpha \\ \bar{Z}_t \text{ is the boundary value of a holomorphic function on domain } D(t) \end{cases} \qquad (2.7)$$

where $D(t)$ is the limit domain of $\Omega^\epsilon(t)$ with boundary $Z = Z(\alpha, t)$, $\alpha \in R$. From the self-similarity law (2.5), we know that at $t = 0$, the velocity $Z_t$ obeys the homogeneity

$$Z_t(\alpha, 0) = \lambda^{s-1} Z_t(\lambda\alpha, 0) \qquad \text{for all } \lambda > 0.$$

Only when $s = 1$, $Z_t$ can be both nontrivial and bounded. For $s = 1$, we see that the gravity term in (2.6) appears negligible. However one could not neglect the surface tension, except where the derivative of the curvature is negligibly small. The role of surface tension will only be clear when we know the profile $Z = Z(\alpha, t)$.

From now on we will focus on finding a profile $Z$, satisfying the gravity free, surface tension free surface water wave equation (2.7) and the self-similarity law (2.5) with $s = 1$. Notice that the self-similarity law for $A$ in this case is

$$A(\alpha, t) = \lambda A(\lambda\alpha, \lambda t) \qquad \text{for all } \lambda > 0. \qquad (2.8)$$

We want the Taylor sign condition hold, that is $A \geq 0$. [1]

## 3. Self-Similar Gravity Free Surface Tension Free Surface Water Waves

Let $Z(\alpha, t) = t\zeta(\beta)$, $Z_t(\alpha, t) = W(\beta)$, $Z_{tt}(\alpha, t) = t^{-1}U(\beta)$ and $A(\alpha, t) = t^{-1}a(\beta)$, where $\beta = \dfrac{\alpha}{t}$ and $a = a(\beta) \geq 0$ for all $\beta \in R$. Assume $(Z, A) = (Z(\alpha, t), A(\alpha, t))$ is a solution of (2.7), then $(\zeta, W, U)$ must satisfy

$$\begin{cases} W(\beta) = \zeta(\beta) - \beta\zeta'(\beta) \\ U(\beta) = -\beta W'(\beta) \\ U(\beta) = ia(\beta)\zeta'(\beta) \\ \overline{W} \text{ is the boundary value of a holomorphic function on } D \end{cases} \qquad (3.1)$$

---

[1] We know this is true if the solution $z = z(\alpha, t)$ satisfies the Taylor sign condition. Using the same argument as in [30], this can also be proved directly for solutions of equation (2.7) satisfying appropriate assumptions at infinity.

where $D$ is the domain on the right as one walks in the direction of increasing $\beta$ on its boundary $\partial D : \zeta = \zeta(\beta), \beta \in R$. Notice that $\overline{W} = \overline{\mathbf{V}} \circ \zeta$, where $\mathbf{V}$ is the profile of velocity field in domain $D$. It is easy to derive from (3.1) that $W' = -\beta \zeta''$, so

$$\beta^2 \zeta''(\beta) = ia(\beta)\zeta'(\beta), \qquad (3.2)$$

and [2]

$$\zeta'(\beta) = e^{i\{b(\beta)+\phi\chi(\beta)\}} \qquad (3.3)$$

where $b = b(\beta)$ is continuous and differentiable on R, and for $\beta \in R$,

$$b'(\beta) = \frac{a(\beta)}{\beta^2}, \quad \chi(\beta) = \begin{cases} 0, & \text{for } \beta < 0 \\ 1, & \text{for } \beta > 0 \end{cases} \qquad (3.4)$$

$\phi$ is a constant[3]. We take $-\pi < \phi \le 0$ for surface water waves. Let $\zeta(0) = 0$. Since $a = a(\beta) \ge 0$, $b = b(\beta)$ is increasing on $R$. Therefore $\zeta = \zeta(\beta), \beta \in R$ is a curve that as $\beta$ goes from $-\infty$ to $0$, concave upwards; at $\beta = 0$ turns downward by an angle of degree $|\phi|$; then continues concave upwardly as $\beta$ goes from $0$ to $\infty$. This appears to be in very good agreement with a common type of wave crests we see in the deeper part of ocean.[4]

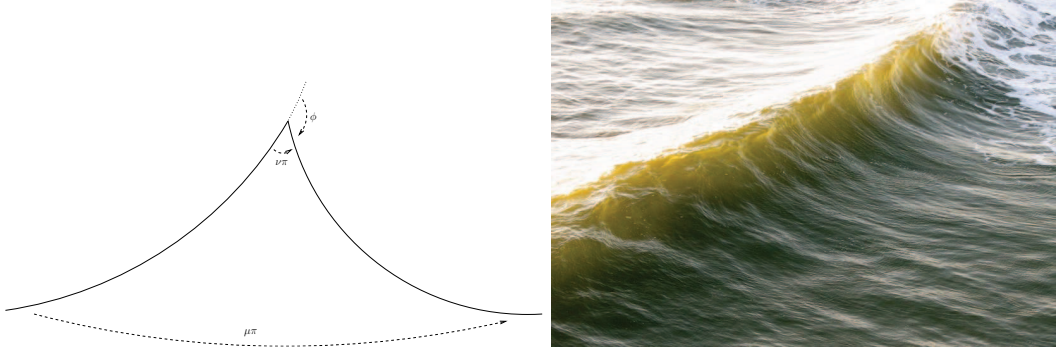

Let

$$\Phi : D \to P_-, \quad \text{with } \Phi(0) = 0$$

be a Riemann mapping from the fluid domain $D$ to the lower half plane $P_-$, and

$$h = h(\beta) := \Phi \circ \zeta = \Phi(\zeta(\beta)) : R \to R.$$

Then $h(0) = 0$; $\zeta \circ h^{-1}(x) = \Phi^{-1}(x)$ and

$$(\zeta \circ h^{-1})' = (h^{-1})' e^{i(b \circ h^{-1} + \phi\chi)} = (\Phi^{-1})'(x), \qquad x \in R \qquad (3.5)$$

---

[2] Without loss of generality we take $\zeta = \zeta(\beta)$ in arclength variable.
[3] It is easy to check $\zeta$ given by (3.3) is a weak solution of (3.2).
[4] The photo is from http://photos.surfline.com.

are the traces on $\partial P_-$ of the holomorphic functions $\Phi^{-1}$ and $\partial_z \Phi^{-1}$ respectively. Now by (3.1), $\overline{W}$ is the trace on $\partial D$ of the holomorphic function $\overline{\mathbf{V}}$ on $D$, that is

$$\overline{W}(\beta) = \overline{\mathbf{V}}(\zeta(\beta)) \quad \beta \in R.$$

We have $\overline{W} \circ h^{-1}(x) = \overline{\mathbf{V}} \circ \Phi^{-1}(x)$, $x \in R$ is the trace on $\partial P_-$ of the holomorphic function $\overline{\mathbf{V}} \circ \Phi^{-1}$, therefore $(\overline{W} \circ h^{-1})'$ and $(\overline{W} \circ h^{-1})'(\zeta \circ h^{-1})'$ are respectively the traces on $\partial P_-$ of the holomorphic functions $\partial_z(\overline{\mathbf{V}} \circ \Phi^{-1})$ and $\partial_z(\overline{\mathbf{V}} \circ \Phi^{-1})\partial_z \Phi^{-1}$. Now $W' = -\beta \zeta'' = -i\beta b'(\beta) e^{i(b+\phi\chi)}$, so

$$(\overline{W} \circ h^{-1})'(\zeta \circ h^{-1})'(x) = i(h^{-1})' h^{-1} (b \circ h^{-1})'(x), \qquad x \in R.$$

Let $0 < \nu\pi \le \pi$ be the angle between the left and right branches of the interface $\zeta = \zeta(\beta), \beta \in R$ about $\beta = 0$; $0 < \mu\pi \le 2\pi$ be the angle between the left and right branches of the interface $\zeta = \zeta(\beta), \beta \in R$ at $\beta = \pm\infty$, so

$$\nu\pi - \phi = \pi, \qquad b(+\infty) - b(-\infty) = (\mu - \nu)\pi. \tag{3.6}$$

The Riemann mapping $\Phi : D \to P_-$ therefore satisfy[5]

$$\Phi^{-1}(z) \sim z^\nu, \quad \text{at } z \sim 0; \qquad \Phi^{-1}(z) \sim z^\mu, \quad \text{at } z \sim \infty.$$

and $(\Phi^{-1})'(z)$ may be assumed to take the form

$$(\Phi^{-1})'(z) = z^{\nu-1}(z-i)^{\mu-\nu} e^{\Psi(z)} \tag{3.7}$$

where $\Psi$ is holomorphic and bounded on $P_-$. Notice that the velocity field $\mathbf{V}$ is harmonic in the fluid domain $D$. Since the fundamental solution of the Laplacian $\frac{1}{2\pi}\ln|z|$ satisfies $|\nabla(\frac{1}{2\pi}\ln|z|)| \lesssim \frac{1}{|z|}$, we assume the velocity field has similar behavior at infinity:

$$|\partial_z \overline{V}(z)| \lesssim |z|^{-1}, \qquad \text{at } z \sim \infty \tag{3.8}$$

We have then

$$|z\partial_z(\overline{\mathbf{V}} \circ \Phi^{-1})\partial_z\Phi^{-1}| = |z\partial_z \overline{\mathbf{V}} \circ \Phi^{-1}(\partial_z\Phi^{-1})^2| \lesssim \left|\frac{z}{\Phi^{-1}}(\partial_z\Phi^{-1})^2\right| \sim |z|^{\mu-1}, \quad \text{at } z \sim \infty$$

On the other hand, $(b \circ h^{-1})'$ is integrable,[6] so we require $|x(\overline{W} \circ h^{-1})'(\zeta \circ h^{-1})'(x)| = |x(h^{-1})' h^{-1}(b \circ h^{-1})'(x)| \lesssim |x|^{2\nu-1}$ in a neighborhood of zero. Now on $\partial P_-$,

$$z\partial_z(\overline{\mathbf{V}} \circ \Phi^{-1})\partial_z\Phi^{-1}\big|_{z=x} = x(\overline{W} \circ h^{-1})'(\zeta \circ h^{-1})'(x) = ix(h^{-1})' h^{-1}(b \circ h^{-1})'(x)$$

is pure imaginary. Therefore $z\partial_z(\overline{\mathbf{V}} \circ \Phi^{-1})\partial_z\Phi^{-1} = i\kappa$ for some constant $\kappa$,[7] and we have

$$(b \circ h^{-1})'(x) = \frac{\kappa}{x(h^{-1})'(x) h^{-1}(x)}. \tag{3.9}$$

---

[5] Here we use the notation $f \sim g$ at $z \sim 0$ to indicate roughly $f/g$ and $g/f$ are bounded in a neighborhood of 0, $f \lesssim g$ means $f \le cg$ for some constant $c$.

[6] We want $(b \circ h^{-1})'$ integrable, since wave phenomenon suggest surface water waves should not role up infinitely near the crest.

[7] We are not considering unusual possibilities.

We know $\kappa > 0$.[8] Notice that
$$h^{-1}(x) \sim x|x|^{\nu-1}, \quad (h^{-1})'(x) \sim |x|^{\nu-1} \quad \text{at } x \sim 0$$
$$h^{-1}(x) \sim x|x|^{\mu-1}, \quad (h^{-1})'(x) \sim |x|^{\mu-1} \quad \text{at } |x| \sim \infty$$

In order for $(b \circ h^{-1})'$ to be integrable on $R$, we must have
$$1/2 < \mu, \quad 0 < \nu < 1/2.$$

The boundary value of the holomorphic function $\Psi$ can be written as $\Psi(x) = g(x) - i\,Hg(x) + c$, $x \in R$, for some real valued function $g$ on $R$, where $H$ is the Hilbert transform:
$$Hf(x) = \frac{1}{\pi} p.v. \int \frac{f(y)}{x-y}\,dy,$$

$c$ is a constant. Hence
$$(\Phi^{-1})'(x) = |x|^{\nu-1}(x^2+1)^{\frac{\mu-\nu}{2}} e^{i(\nu-1)\pi(\chi(x)-1)} e^{i(\mu-\nu)\arctan\frac{-1}{x}} e^{g(x) - i\,Hg(x) + c} \quad (3.10)$$

Sum up (3.9), (3.5), (3.10), (3.6), we arrive at the following equation for $\zeta'(x) = e^{i(b(x)+\phi\chi(x))}$:

Let
$$1/2 < \mu \leq 2, \quad 0 < \nu < 1/2. \quad (3.11)$$

For some real valued function $g$ on $R$,
$$\begin{cases} (b \circ h^{-1})'(x) = \frac{\kappa}{xh^{-1}(x)(h^{-1})'(x)} \\ (h^{-1})'(x) = |x|^{\nu-1}(x^2+1)^{\frac{\mu-\nu}{2}} e^{g(x)} \\ g = H(b \circ h^{-1}(x) - (\mu-\nu)\arctan\frac{-1}{x}) \end{cases} \quad (3.12)$$

where $h^{-1}(0) = 0$, $\kappa > 0$ is the constant so that
$$\int_{-\infty}^{\infty} \frac{\kappa}{xh^{-1}(x)(h^{-1})'(x)}\,dx = (\mu-\nu)\pi.$$

We prove the following results in this paper. We have

**Theorem 3.1** (Apriori Estimate). *Let $1/2 < \mu \leq 2$, $0 < \nu < 1/2$. Let $G(x) = g(x) + \frac{1}{2}(\mu-\nu)\ln(x^2+1)$. Assume that $G$ is even, increasing on $[0, \infty)$, and $g \in C(R) \cap L^\infty(R)$ is a solution of the system (3.12). Then there exist $x_0 > 0$, and a constant $c(\mu, \nu)$, depending only on $\mu$, $\nu$, such that*
$$\|G(x) - \frac{1}{2}(\mu-\nu)\ln(x^2 + x_0^2)\|_{L^\infty(R)} \leq c(\mu, \nu)$$

**Remark 3.2.** Basically, Theorem 3.1 states that the Riemann Mapping $\Phi$ has the property that
$$(\Phi^{-1})'(z) = z^{\nu-1}(z - x_0 i)^{\mu-\nu} e^{\Psi(z)}$$

---

[8]This is because $(b \circ h^{-1})'(x) \geq 0$.

where $x_0 > 0$ and $\Psi$ is bounded by a universal constant $c(\mu, \nu)$. Notice that for any $\lambda > 0$, the rescaled function $\Phi_\lambda^{-1} = \Phi_\lambda^{-1}(z) := \Phi^{-1}(\lambda z) : z \in P_- \to D$ is also a Riemann Mapping. From the proof, we will see that $x_0 > 0$ is the number satisfying

$$b \circ h^{-1}(x_0) - b \circ h^{-1}(-x_0) = (\frac{1}{2} - \nu)\pi.$$

**Theorem 3.3** (Existence)**.** *Let $1/2 < \mu \leq 2$, $0 < \nu < 1/2$. In $C(R) \cap L^\infty(R)$ there exists a solution $g$ of the system (3.12), such that the function $G = G(x) = g(x) + \frac{1}{2}(\mu - \nu)\ln(x^2 + 1) : R \to R$ is even and increasing on $[0, \infty)$. Moreover, $xg' \in L^\infty(R)$.*

**Remark 3.4.** Given $g \in L^\infty(R) \cap C(R)$ a solution of the system (3.12), we can construct a solution for system (3.1) and then the systems (2.7) and (1.1)(with $\mathfrak{g} = 0$) through the following procedure. First we construct from (3.12) the homeomorphism $h^{-1} : R \to R$ with $h^{-1}(0) = 0$ and the increasing function $b$, satisfying $b(\infty) = (\mu - \nu)\arctan\frac{-1}{x}\big|_{x=\infty}$.[9] We then obtain the curve $\zeta = \zeta(\beta)$ satisfying $\zeta(0) = 0$ and $\zeta'(\beta) = e^{i(b(\beta) + (\nu-1)\pi\chi(\beta))}$, for $\beta \in R$. Let $a(\beta) = \beta^2 b'(\beta)$. We know $a(\beta) \geq 0$. Let $W = \zeta(\beta) - \beta\zeta'(\beta)$, $U = -\beta W'(\beta)$. It is easy to check $(\zeta, W, U, a)$ satisfies the first three equations in system (3.1). Now let $\Psi$ be the bounded holomorphic function on $P_-$ taking boundary value $g - iHg + (\nu - 1)\pi$, and let $\Xi$ be the holomorphic function on $P_-$, with $\Xi(0) = 0$, and

$$\Xi'(z) = z^{\nu-1}(z - i)^{\mu-\nu} e^{\Psi(z)}.$$

Notice that $\Xi'(z) \neq 0$ for $z \in P_-$, therefore $\Xi$ is conformal on $P_-$. It is easy to check $\Xi'(x) = (\zeta \circ h^{-1})'(x)$. Therefore $\Xi(x) = \zeta \circ h^{-1}(x)$ for $x \in R$. It then follows from the argument principle that $\Xi$ is a Riemann Mapping from $P_-$ to the domain $D$. Now $W' = -\beta\zeta''(\beta) = -i\beta b'(\beta)e^{i(b+(\nu-1)\pi\chi)}$ therefore

$$(\overline{W} \circ h^{-1})'(x) = \frac{i\kappa}{x\Xi'(x)} \qquad x \in R$$

so $(\overline{W} \circ h^{-1})'$ is the trace on $\partial P_-$ of the holomorphic function $\frac{i\kappa}{z\Xi'(z)}$, consequently $\overline{W} \circ h^{-1}$ is the trace on $\partial P_-$ of the antiderivative of $\frac{i\kappa}{z\Xi'(z)}$, we name it $\Lambda$, which is also holomorphic on $P_-$. $\overline{W}$ is then the boundary value of the holomorphic function $\Lambda \circ \Xi^{-1}$ on $D$. This shows that $(\zeta, W, U, a)$ constructed as above is a solution for the surface water wave system (3.1). $(Z, A) = (t\zeta(\frac{\alpha}{t}), t^{-1}a(\frac{\alpha}{t}))$, $\alpha \in R$, $t \neq 0$ is then a solution for (2.7). From the equivalence of (2.7) and (1.1)(with $\mathfrak{g} = 0$), we have a self-similar solution for (2.7) and (1.1)(with $\mathfrak{g} = 0$). Moreover this solution satisfies the Taylor sign condition. That is $a \geq 0$ or $-\frac{\partial P}{\partial \mathbf{n}} \geq 0$ on the interface.

---

[9] We can set $b(\infty)$ to equal to any number. The difference a different $b(\infty)$ makes is a rotation of the interface.

## 4. Discussion

Before we prove Theorems 3.1 and 3.3, we give a brief discussion on the self-similar gravity free surface tension free surface water wave we found in section 3. First, we want to understand the effect of the surface tension to this wave. Let $Z = Z(\alpha, t) = t\zeta(\frac{\alpha}{t})$, $\alpha \in R$, where $\zeta$ is the wave profile found in Theorem 3.3. Let $\beta = \frac{\alpha}{t}$. We have $Z_\alpha = \zeta'(\beta)$, $Z_{\alpha\alpha} = t^{-1} i b'(\beta)\zeta'(\beta)$, for $\beta \neq 0$, and we calculate and find the surface tension term for $Z$ in equation (2.2) is

$$-\sigma \frac{Z_\alpha}{|Z_\alpha|^2} \partial_\alpha (\operatorname{Im} \frac{\bar{Z}_\alpha Z_{\alpha\alpha}}{|Z_\alpha|^3}) = -\sigma \zeta'(\beta) t^{-2} b''(\beta), \qquad \text{for } \beta \neq 0.$$

Now we calculate $b''$ from equation (3.12). We use the following notations. We write $f \simeq g$ if $1/c \leq f/g \leq c$ for some positive constant $c$; we write $f \lesssim g$ if $f \leq cg$ for some positive constant $c$. First we mention that from (5.4) in section 5, we have

$$h^{-1}(x) \simeq x(h^{-1})'(x), \qquad \text{for } x \in R \tag{4.1}$$

We know

$$b' \circ h^{-1} = \frac{(b \circ h^{-1})'}{(h^{-1})'} = \frac{\kappa}{xh^{-1}(h^{-1})'^2},$$

therefore

$$b' \circ h^{-1}(x) \simeq \begin{cases} |x|^{1-3\nu}, & \text{for } |x| \leq 1 \\ |x|^{1-3\mu}, & \text{for } |x| \geq 1. \end{cases}$$

Also

$$|b'' \circ h^{-1}| = \left|\frac{1}{(h^{-1})'}\left(\frac{(b \circ h^{-1})'}{(h^{-1})'}\right)'\right| = \left|\frac{(b \circ h^{-1})'}{(h^{-1})'^2}\left(-\frac{1}{x} - \frac{(h^{-1})'}{h^{-1}} - 2\frac{(h^{-1})''}{(h^{-1})'}\right)\right| \lesssim \frac{(b \circ h^{-1})'}{|x|(h^{-1})'^2}$$

then

$$|b'' \circ h^{-1}(x)| \lesssim \frac{1}{x^2 |h^{-1}(x)|(h^{-1})'^3} \lesssim \begin{cases} |x|^{1-4\nu} & \text{for } |x| \leq 1 \\ |x|^{1-4\mu} & \text{for } |x| \geq 1. \end{cases}$$

Now from (4.1) and Theorem 3.3,

$$x = h(\beta) \simeq \begin{cases} \beta|\beta|^{\frac{1}{\nu}-1} & \text{for } |\beta| \leq 1 \\ \beta|\beta|^{\frac{1}{\mu}-1} & \text{for } |\beta| \geq 1 \end{cases}$$

therefore

$$b'(\beta) \simeq \begin{cases} |\beta|^{\frac{1}{\nu}-3} & \text{for } |\beta| \leq 1 \\ |\beta|^{\frac{1}{\mu}-3} & \text{for } |\beta| \geq 1 \end{cases}$$

$$|b''(\beta)| \lesssim \begin{cases} |\beta|^{\frac{1}{\nu}-4} & \text{for } |\beta| \leq 1 \\ |\beta|^{\frac{1}{\mu}-4} & \text{for } |\beta| \geq 1 \end{cases}$$

We see the strength of the surface tension term is

$$\left| -\sigma \frac{Z_\alpha}{|Z_\alpha|^2} \partial_\alpha (\operatorname{Im} \frac{\bar{Z}_\alpha Z_{\alpha\alpha}}{|Z_\alpha|^3}) \right| \lesssim \begin{cases} t^{2-\frac{1}{\nu}} |\alpha|^{\frac{1}{\nu}-4} & \text{for } 0 < |\alpha| \leq t \\ t^{2-\frac{1}{\mu}} |\alpha|^{\frac{1}{\mu}-4} & \text{for } |\alpha| \geq t \end{cases}$$

and the acceleration

$$|Z_{tt}| = |t^{-1}U(\beta)| = |t^{-1}\beta^2 b'(\beta)| \simeq \begin{cases} t^{-\frac{1}{\nu}}|\alpha|^{\frac{1}{\nu}-1} & \text{for } |\alpha| \leq t \\ t^{-\frac{1}{\mu}}|\alpha|^{\frac{1}{\mu}-1} & \text{for } |\alpha| \geq t \end{cases}$$

Comparing with the acceleration $Z_{tt}$ and the term $iAZ_\alpha = t^{-1}ia\zeta'(\frac{\alpha}{t})$, for $t$ close enough to zero, the effect from surface tension is negligibly small where $|\alpha| \gtrsim t^{2/3-}$. The effect of surface tension is significant where $|\alpha| \lesssim t^{2/3}$.

Nevertheless, this discussion and the heuristic one in section 2 deserve to be made rigorous. We plan to study the asymptotic stability of the surface water wave found in Theorem 3.3 and the effect of gravity and surface tension in rigorous mathematical terms in upcoming works.

We can also calculate the velocity profile $W$ for the solution found in Theorem 3.3. First we take $b(0) = -\frac{\phi}{2} = -\frac{\nu-1}{2}\pi$, so that the phase $b + \phi\chi$ of $\zeta'$ is an odd function, and the wave profile $\zeta = \zeta(\beta)$, $\beta \in R$ is symmetric about the vertical axis $\beta = 0$. It is easy to see from $W' = -i\beta b'(\beta)e^{i(b+\phi\chi)}$ that the horizontal velocity $\operatorname{Re} W$ is an odd function, and the vertical velocity $\operatorname{Im} W$ is even. We have

$$\beta b'(\beta) \simeq \begin{cases} \beta|\beta|^{\frac{1}{\nu}-3} & \text{for } |\beta| \leq 1 \\ \beta|\beta|^{\frac{1}{\mu}-3} & \text{for } |\beta| \geq 1 \end{cases}$$

and

$$b(\infty) - b(\beta) \simeq \beta^{\frac{1}{\mu}-2}, \qquad \text{for } \beta \geq 1.$$

So for $\beta$ in any bounded interval, $W$ is bounded. Notice that $b(\infty) + \phi = \frac{\mu-1}{2}\pi$. For $\beta > 1$ We calculate

$$W(\beta) - W(1) = -ie^{i(\frac{\mu-1}{2}\pi)} \int_1^\beta \gamma b'(\gamma) e^{i(b(\gamma)-b(\infty))} \, d\gamma$$

and find that as $\beta \to \infty$, $W(\beta)$ asymptotically points in the direction $-ie^{i(\frac{\mu-1}{2}\pi)}$ and has magnitude $O(\beta^{\frac{1}{\mu}-1})$ if $\mu < 1$. When $\mu = 1$, the horizontal velocity $\operatorname{Re} W$ is bounded, the vertical velocity points in the direction $-i$, and has magnitude $O(\ln \beta)$. For $\mu > 1$, $W$ is bounded for all $\beta \in R$.[10]

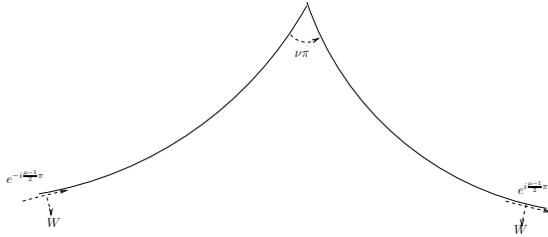
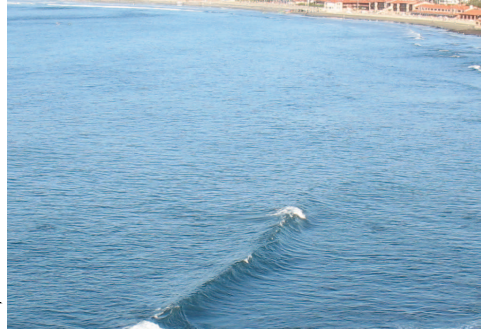

---

[10]This photo is taken by the author at San Diego, CA.

## 5. The Proof Of The Main Theorems

We prove Theorem 3.3 by Schaefer's fixed point theorem. We use the following version, which is a slight modification from the one stated in Chapter 9 of [13]. The same proof as in [13] works for this one as well.

**Theorem 5.1** (Schaefer's Fixed Point Theorem). *Let $X$ be a Banach space, $X_1$ be a convex subset of $X$, and $0 \in X_1$. Suppose $T : X \to X$ is continuous and compact, and $T : X_1 \to X_1$. Assume further that the set*

$$\{u \in X_1 \,|\, u = \lambda T[u] \text{ for some } 0 \leq \lambda \leq \lambda_0\}$$

*is bounded. Here $\lambda_0 > 0$ is a constant. Then $T$ has a fixed point in $X_1$.*

Let $\mu > 1/2$, $0 < \nu < 1/2$. Let

$$X = \{g \in C(R) \cap L^\infty(R) \,|\, \lim_{x \to \pm\infty} g(x) = 0\},$$

with norm $\|g\| = \|g\|_{L^\infty(R)} = \|g\|_\infty$. For $g \in X$, we define $T[g]$ by:

$$T[g] = H(F - (\mu - \nu) \arctan \frac{-1}{x}) \tag{5.1}$$

where $H$ is the Hilbert transform, $h^{-1}(0) = 0$,

$$\begin{cases} F'(x) = \frac{\kappa}{xh^{-1}(x)(h^{-1})'(x)} \\ (h^{-1})'(x) = |x|^{\nu-1}(x^2+1)^{\frac{\mu-\nu}{2}}e^{g(x)} \\ F(\infty) = (\mu - \nu) \arctan \frac{-1}{x}\big|_{x=\infty} \\ F(-\infty) = (\mu - \nu) \arctan \frac{-1}{x}\big|_{x=-\infty} \end{cases} \tag{5.2}$$

Notice that $H1 = 0$. For $x \in R$, let $G(x) = g(x) + \frac{1}{2}(\mu - \nu)\ln(x^2 + 1)$,

$$X_1 = \{g \in X \,|\, G \text{ is even},\ 0 \leq G(x) - G(y) \leq (\mu - \nu)(\ln \frac{x}{y} + 2),\ \text{for all } 0 \leq y \leq x\}$$

For $I \subset R$ an interval, let $\fint_I g = \frac{1}{|I|} \int_I g$. In what follows we use the notations specified above without further clarifying. $c(\mu, \nu)$, $c(\mu, \nu, M)$ etc. are constants depending on $\mu, \nu$ or $\mu, \nu, M$, they need not be the same in different contexts.

### 5.1. The continuity and compactness of $T$. We have the following

**Lemma 5.2.** *1. Let $g \in X$. We have for all $x \in R$,*

$$|h^{-1}(x)| \leq \frac{1}{\nu}|x|^\nu(x^2+1)^{\frac{1}{2}(\mu-\nu)}e^{\|g\|_\infty}$$
$$|h^{-1}(x)| \geq e^{-(\mu-1)}|x|^\nu(x^2+1)^{\frac{1}{2}(\mu-\nu)}e^{\fint_0^x g} \geq e^{-(\mu-1)}|x|^\nu(x^2+1)^{\frac{1}{2}(\mu-\nu)}e^{-\|g\|_\infty} \tag{5.3}$$

*2. Let $g \in X_1$. We have for $x \in R$,*

$$\frac{1}{x}h^{-1}(x) \leq \frac{1}{\nu}|x|^{\nu-1}(x^2+1)^{\frac{1}{2}(\mu-\nu)}e^{g(x)}$$
$$\frac{1}{x}h^{-1}(x) \geq e^{1-\nu-6(\mu-\nu)}|x|^{\nu-1}(x^2+1)^{\frac{1}{2}(\mu-\nu)}e^{g(x)} \tag{5.4}$$

*Proof.* The first inequalities in (5.3) and (5.4) are straightforward. We prove the second inequalities using Jensen's inequality. Let $x > 0$. We have

$$h^{-1}(x) = \int_0^x (h^{-1})'(y)\, dy \geq x e^{\fint_0^x ((\nu-1)\ln|y| + \frac{1}{2}(\mu-\nu)\ln(y^2+1) + g(y))\, dy} = x e^{\fint_0^x ((\nu-1)\ln|y| + G(y))\, dy}$$

Using integration by parts, we get

$$\int_0^x ((\nu-1)\ln|y| + \frac{1}{2}(\mu-\nu)\ln(y^2+1))\, dy = x\{(\nu-1)\ln|x| + \frac{1}{2}(\mu-\nu)\ln(x^2+1)\}$$

$$- \int_0^x (\nu - 1 + (\mu-\nu)\frac{y^2}{y^2+1})\, dy \geq x\{(\nu-1)\ln|x| + \frac{1}{2}(\mu-\nu)\ln(x^2+1)\} - (\mu-1)x$$

The second inequality in (5.3) therefore follows. For the second inequality in (5.4) and $g \in X_1$ we estimate

$$0 \leq \fint_{x/2}^x G - \fint_0^x G = \frac{1}{x}\left(\int_{x/2}^x G(y_1)\, dy_1 - \int_0^{x/2} G(y_2)\, dy_2\right)$$

$$\leq (\mu - \nu)\{\frac{1}{x}\left(\int_{x/2}^x \ln y_1\, dy_1 - \int_0^{x/2} \ln y_2\, dy_2\right) + 1\}$$

$$\leq (\mu - \nu)(\ln 2 + 1)$$

therefore

$$\fint_0^x G \geq \fint_{x/2}^x G - (\mu - \nu)(\ln 2 + 1)$$

$$\geq G(x/2) - (\mu-\nu)(\ln 2 + 1) \geq G(x) - (\mu-\nu)(2\ln 2 + 3)$$

The second inequality in (5.4) therefore follows. $\square$

**Lemma 5.3.** *Let $g \in X$. There exists a constant $c(\mu,\nu) > 0$, depending on $\mu,\nu$, such that*

$$\frac{1}{c(\mu,\nu)} e^{-2\|g\|_\infty} \leq \kappa \leq c(\mu,\nu) e^{2\|g\|_\infty}$$

This is straightforward, we omit the proof.

**Lemma 5.4.** *Let $g \in X$, $f(x) = F(x) - (\mu-\nu)\arctan\frac{-1}{x}$. We have*

1.
$$|f(x)| \leq c(\mu,\nu) e^{4\|g\|_\infty}|x|^{1-2\mu} + (\mu-\nu)|x|^{-1}, \quad \text{for } |x| \geq 1$$
$$|f(x)| \leq 2(\mu-\nu)\pi, \quad \text{for all } x \in R \quad (5.5)$$
$$|f'(x)| \leq c(\mu,\nu) e^{4\|g\|_\infty}\frac{1}{|x|^{2\nu}(x^2+1)^{\mu-\nu}} + (\mu-\nu)\frac{1}{x^2+1}, \quad \text{for all } x$$

*where $c(\mu,\nu)$ is a constant depending on $\mu,\nu$.*

2. *Assume $\|g\|_\infty \leq M$. Then*

$$|T[g](x)| \leq c(\mu,\nu,M)(|x|^{1/2-\mu} + |x|^{-1/2}), \quad \text{for } |x| \geq 4 \quad (5.6)$$

*where $c(\mu,\nu,M)$ is a constant depending on $\mu,\nu,M$.*

*Proof.* (5.5) is straightforward, we omit the proof.

We use (5.5) to prove (5.6). Assume $x \geq 4$. Notice that

$$T[g](x) = Hf(x) = \frac{1}{\pi}\int_{|x-y|\leq 1} \frac{f(y) - f(x)}{x-y}\,dy + \frac{1}{\pi}\int_{|x-y|>1} \frac{f(y)}{x-y}\,dy = I + II$$

we have

$$|I| \leq \frac{1}{\pi}\int_{|x-y|\leq 1} \frac{1}{|x-y|}|\int_x^y f'(\gamma)\,d\gamma|\,dy \leq c(\mu,\nu,M)(|x|^{-2\mu} + |x|^{-2})$$

We further decompose

$$II = \frac{1}{\pi}\int_{1<|x-y|\leq \frac{x}{2}} \frac{f(y)}{x-y}\,dy + \frac{1}{\pi}\int_{|x-y|>\frac{x}{2}} \frac{f(y)}{x-y}\,dy = II_1 + II_2$$

and we know $|y| \geq |x| - |x-y|$, so

$$|II_1| \leq c(\mu,\nu,M)(|x|^{1-2\mu} + |x|^{-1})\ln|x|$$

Now

$$II_2 = \frac{1}{\pi}\int_{y>\frac{3}{2}x, y<-\frac{1}{2}x} \frac{f(y)}{x-y}\,dy + \frac{1}{\pi}\int_{-\frac{1}{2}x<y<\frac{1}{2}x} \frac{f(y)}{x-y}\,dy = II_{21} + II_{22}$$

where since for $y > \frac{3}{2}x$, $y - x > \frac{1}{4}y > 0$, and for $y < -\frac{1}{2}x$, we have $y - x < y < 0$, therefore

$$|II_{21}| \leq c(\mu,\nu,M)(|x|^{1-2\mu} + |x|^{-1})$$

For $II_{22}$, we further decompose (notice that $x^{1/2} \leq \frac{1}{2}x$)

$$|II_{22}| = |\frac{1}{\pi}\int_{|y|\leq x^{1/2}} \frac{f(y)}{x-y}\,dy + \frac{1}{\pi}\int_{x^{1/2}\leq |y|<\frac{1}{2}x} \frac{f(y)}{x-y}\,dy|$$
$$\leq C(\mu,\nu,M)(|x|^{-1/2} + |x|^{1/2-\mu})$$

This proves (5.6) for $x \geq 4$. The proof for $x \leq -4$ is similar, we omit. $\square$

The above three Lemmas is sufficient for the proof of the compactness of $T: X \to X$. We need the following Lemma for the continuity of $T$.

For $i = 1,2$, let $g_i \in X$, $(h_i^{-1})'(x) = |x|^{\nu-1}(x^2+1)^{\frac{1}{2}(\mu-\nu)}e^{g_i(x)}$,

$$F_i'(x) = \frac{\kappa_i}{xh_i^{-1}(x)(h_i^{-1})'(x)}$$

with $h_i^{-1}(0) = 0$, $F_i(\infty) - F_i(-\infty) = (\mu-\nu)\pi$, $F_i(\infty) = (\mu-\nu)\arctan\frac{-1}{x}\big|_{x=\infty}$, and $\|g_i\|_\infty \leq M$.

**Lemma 5.5.** *We have*

$$|F_1'(x) - F_2'(x)| \leq c(\mu,\nu,M)\|g_2 - g_1\|_\infty \frac{1}{|x|^{2\nu}(x^2+1)^{\mu-\nu}}, \quad \text{for } x \in R$$

$$|F_1(x) - F_2(x)| \leq c(\mu,\nu,M)\|g_2 - g_1\|_\infty \frac{1}{|x|^{1-2\mu}}, \quad \text{for } |x| \geq 1 \qquad (5.7)$$

$$|F_1(x) - F_2(x)| \leq c(\mu,\nu,M)\|g_2 - g_1\|_\infty, \quad \text{for } x \in R.$$

*Proof.* We know
$$F_1'(x) - F_2'(x) = \frac{\kappa_1 - \kappa_2}{xh_1^{-1}(x)(h_1^{-1})'(x)} + \frac{\kappa_2}{xh_2^{-1}(x)(h_2^{-1})'(x)}\left(\frac{h_2^{-1}(x)(h_2^{-1})'(x)}{h_1^{-1}(x)(h_1^{-1})'(x)} - 1\right)$$

Now
$$\frac{(h_2^{-1})'(x)}{(h_1^{-1})'(x)} = e^{g_2(x) - g_1(x)}$$
$$\frac{(h_2^{-1})(x)}{(h_1^{-1})(x)} = \frac{(h_2^{-1})'(\theta)}{(h_1^{-1})'(\theta)} = e^{g_2(\theta) - g_1(\theta)}$$

for some $0 < |\theta| < |x|$. Here we used the mean value Theorem. So
$$\left|\frac{h_2^{-1}(x)(h_2^{-1})'(x)}{h_1^{-1}(x)(h_1^{-1})'(x)} - 1\right| \leq 2e^{2\|g_2 - g_1\|_\infty} \|g_2 - g_1\|_\infty$$

Therefore from $\int_{-\infty}^{\infty}(F_1'(x) - F_2'(x))\,dx = 0$ we have
$$\frac{|\kappa_1 - \kappa_2|}{\kappa_1} \leq 2e^{2\|g_2 - g_1\|_\infty} \|g_2 - g_1\|_\infty$$

and
$$|F_1'(x) - F_2'(x)| \leq c(\mu, \nu, M) \|g_2 - g_1\|_\infty \frac{1}{|x|^{2\nu}(x^2 + 1)^{\mu - \nu}}$$

for some constant $c(\mu, \nu, M)$ depending on $\mu, \nu, M$. The second and last inequalities in (5.7) follow by integration. □

We are now ready to prove

**Proposition 5.6.** $T : X \to X$ *is compact and continuous.*

*Proof.* Recall that for $1 < p < \infty$, the Hilbert transform $H : L^p(R) \to L^p(R)$ is bounded. Let $g \in X$, $f = F - (\mu - \nu)\arctan\frac{-1}{x}$ be as in Lemma 5.4. We know $T[g] = H(f)$. Take $1 < r, q < \infty$ so that $(2\mu - 1)r > 1$ and $2\nu q < 1$. From Lemma 5.4, we have $f \in L^r(R)$, and $f' \in L^q(R)$. Define
$$Y(I) = \{f \in L^r(I) \mid f' \in L^q(I)\}$$

with norm $\|f\|_{Y(I)} = \|f\|_{L^r(I)} + \|f'\|_{L^q(I)}$. Therefore for $g \in X$, $T[g] \in Y(R)$, and
$$\|T[g]\|_{Y(R)} \leq c(\mu, \nu, \|g\|_\infty),$$

and from Lemma 5.5, for $g_1, g_2 \in X$, with $\|g_i\|_\infty \leq M$,
$$\|T[g_1] - T[g_2]\|_{Y(R)} \leq c(\mu, \nu, M)\|g_1 - g_2\|_\infty.$$

Now it follows from the same proof as for Sobolev embeddings (see [13]) that $Y(R) \subset X$ and the embedding is continuous. Therefore $T : X \to X$ is continuous. The same proof for Sobolev compact embeddings (i.e. use Arzela-Ascoli Theorem, see [13]) also shows that $Y(I)$ is compactly embedded in $C(I)$ for any bounded interval $I \subset R$. Using (5.6) and a standard diagonal argument, we have that $T : X \to X$ is compact.

□

5.2. **A priori estimate.** First we have the following Lemma

**Lemma 5.7.** *Let $\mathfrak{F} \in C^1(R)$. Assume $\mathfrak{F}'$ is even, nonnegative, and decreasing on $[0, \infty)$; and $\mathfrak{F}(\infty) - \mathfrak{F}(-\infty) = (\mu - \nu)\pi$. Then the function $\mathfrak{G} = H\mathfrak{F}(x) - H\mathfrak{F}(0)$ is even, and is increasing on $[0, \infty)$. Moreover for $y > x \geq 0$,*

$$0 \leq H\mathfrak{F}(y) - H\mathfrak{F}(x) = \frac{1}{\pi} \int_x^y \frac{\mathfrak{F}(\gamma) - \mathfrak{F}(-\gamma)}{\gamma} d\gamma + \mathfrak{R}(x, y) \tag{5.8}$$

*where*

$$-2(\mu - \nu) \leq \mathfrak{R}(x, y) \leq 2(\mu - \nu)$$

*Proof.* From the assumptions it is easy to show that the function $\mathfrak{G} = H\mathfrak{F}(x) - H\mathfrak{F}(0)$ is even, and is increasing on $[0, \infty)$. For $y > x \geq 0$, we know

$$H\mathfrak{F}(y) - H\mathfrak{F}(x) = \frac{1}{\pi} \int_0^\infty \frac{\mathfrak{F}(y - \gamma) - \mathfrak{F}(y + \gamma) - \mathfrak{F}(x - \gamma) + \mathfrak{F}(x + \gamma)}{\gamma} d\gamma$$

$$= \frac{1}{\pi} \int_0^y \frac{\mathfrak{F}(y - \gamma) - \mathfrak{F}(y + \gamma)}{\gamma} d\gamma - \frac{1}{\pi} \int_0^x \frac{\mathfrak{F}(x - \gamma) - \mathfrak{F}(x + \gamma)}{\gamma} d\gamma$$

$$- \frac{1}{\pi} \int_x^y \frac{\mathfrak{F}(x - \gamma) - \mathfrak{F}(x + \gamma)}{\gamma} d\gamma$$

$$+ \frac{1}{\pi} \int_y^\infty \frac{\mathfrak{F}(y - \gamma) - \mathfrak{F}(y + \gamma) - \mathfrak{F}(x - \gamma) + \mathfrak{F}(x + \gamma)}{\gamma} d\gamma$$

Now

$$\int_y^\infty \frac{\mathfrak{F}(y - \gamma) - \mathfrak{F}(y + \gamma) - \mathfrak{F}(x - \gamma) + \mathfrak{F}(x + \gamma)}{\gamma} d\gamma$$

$$= \int_y^\infty \frac{1}{\gamma} \left( \int_x^y \mathfrak{F}'(\tau - \gamma) d\tau - \int_x^y \mathfrak{F}'(\tau + \gamma) d\tau \right) d\gamma$$

and

$$0 \leq \int_y^\infty \frac{1}{\gamma} \int_x^y \mathfrak{F}'(\tau - \gamma) d\tau \, d\gamma \leq \frac{1}{y} \int_x^y \int_y^\infty \mathfrak{F}'(\tau - \gamma) d\gamma \, d\tau \leq \frac{\mu - \nu}{2}\pi$$

similarly

$$0 \leq \int_y^\infty \frac{1}{\gamma} \int_x^y \mathfrak{F}'(\tau + \gamma) d\tau \, d\gamma \leq \frac{\mu - \nu}{2}\pi$$

Also,

$$0 \leq \int_0^y \frac{\mathfrak{F}(y + \gamma) - \mathfrak{F}(y - \gamma)}{2\gamma} d\gamma = \int_0^y \left( \fint_{y-\gamma}^{y+\gamma} \mathfrak{F}' \right) d\gamma \leq \int_0^y \mathfrak{F}'(y - \gamma) d\gamma \leq \frac{\mu - \nu}{2}\pi$$

Now we further rewrite

$$\frac{1}{\pi} \int_x^y \frac{\mathfrak{F}(x - \gamma) - \mathfrak{F}(x + \gamma)}{\gamma} d\gamma = \frac{1}{\pi} \int_x^y \frac{\mathfrak{F}(-\gamma) - \mathfrak{F}(\gamma)}{\gamma} d\gamma$$

$$+ \frac{1}{\pi} \int_x^y \frac{\mathfrak{F}(x - \gamma) - \mathfrak{F}(-\gamma) - \mathfrak{F}(x + \gamma) + \mathfrak{F}(\gamma)}{\gamma} d\gamma$$

and we have

$$0 \leq \int_x^y \frac{\mathfrak{F}(x - \gamma) - \mathfrak{F}(-\gamma)}{\gamma} d\gamma = \int_x^y \frac{1}{\gamma} \int_0^x \mathfrak{F}'(\tau - \gamma) d\tau \, d\gamma$$

$$\leq \frac{1}{x} \int_0^x (\mathfrak{F}(\tau - x) - \mathfrak{F}(\tau - y)) d\tau \leq \frac{\mu - \nu}{2}\pi$$

similarly
$$0 \le \int_x^y \frac{\mathfrak{F}(x+\gamma) - \mathfrak{F}(\gamma)}{\gamma} d\gamma \le \frac{\mu - \nu}{2}\pi$$
Sum up the above calculation we have (5.8)

□

Notice that $H(\arctan \frac{-1}{x}) - H(\arctan \frac{-1}{x})(0) = \frac{1}{2}\ln(x^2+1)$. Therefore
$$T[g](x) - T[g](0) + \frac{1}{2}(\mu - \nu)\ln(x^2+1) = HF(x) - HF(0)$$
and from Lemma 5.7, we have
$$T : X_1 \to X_1.$$

We use the following notations: we use $\mathfrak{R}$ to indicate a function that satisfies
$$c_1(\mu, \nu) \le \mathfrak{R} \le c_2(\mu, \nu)$$
for some constants $c_i(\mu, \nu)$ depending only on $\mu, \nu$. $e^{\mathfrak{R}}$ indicates a function that is bounded above and below by two positive constants depending only on $\mu, \nu$. We often do the following calculation: Assume $c_1(\mu, \nu) \le \mathfrak{R} \le c_2(\mu, \nu)$. For a function $f \in C^1(I)$, satisfying $f' \ge 0$ on interval $I$, we know for $x_1 < x_2$, $x_1, x_2 \in I$,
$$c_1(\mu,\nu)(f(x_2) - f(x_1)) \le \int_{x_1}^{x_2} e^{\mathfrak{R}(x)} f'(x)\, dx \le c_2(\mu,\nu)(f(x_2) - f(x_1)),$$
we therefore simply write
$$\int_{x_1}^{x_2} e^{\mathfrak{R}} f'(x)\, dx = e^{\mathfrak{R}}(f(x_2) - (f(x_1))).$$
The two $\mathfrak{R}$'s are not necessarily the same. In general, $\mathfrak{R}$'s and $e^{\mathfrak{R}}$'s appearing in different contexts need not be the same.

We have the following a priori estimate for $T : X_1 \to X_1$.

**Proposition 5.8.** *Suppose $g \in X_1$ satisfies*
$$g = \lambda T[g], \qquad \text{for some } 0 < \lambda \le \min\{\frac{1-2\nu}{2(\mu-\nu)}, \frac{2\mu-1}{2(\mu-\nu)}\}. \tag{5.9}$$
*Then $\|g\|_\infty \le c(\mu, \nu)$, where $c(\mu, \nu)$ is a constant depending only on $\mu, \nu$.*

*Proof.* Let $x > 0$, and $0 < \lambda \le \min\{\frac{1-2\nu}{2(\mu-\nu)}, \frac{2\mu-1}{2(\mu-\nu)}\}$. Let $g \in X_1$ and $F$ be defined by (5.2). Notice that $F$ satisfies the assumption of Lemma 5.7. Therefore
$$HF(x) - HF(0) = \frac{1}{\pi} \int_0^x \frac{F(y) - F(-y)}{y} dy + \mathfrak{R}.$$
Let
$$u(x) = \frac{1}{\pi} \int_0^x \frac{F(y) - F(-y)}{y} dy \tag{5.10}$$
Then
$$T[g](x) = T[g](0) - \frac{1}{2}(\mu - \nu)\ln(x^2+1) + u(x) + \mathfrak{R} \tag{5.11}$$

Assume $g$ satisfies (5.9). Then
$$g(x) = \lambda(T[g](0) - \frac{1}{2}(\mu - \nu)\ln(x^2 + 1) + u(x) + \mathfrak{R}) \tag{5.12}$$

Now from (5.10), we have
$$xu'(x) = \frac{1}{\pi}(F(x) - F(-x)), \qquad (xu'(x))' = \frac{2}{\pi}F'(x) \tag{5.13}$$

where by (5.4)
$$F'(x) = \frac{\kappa e^{\mathfrak{R}}}{e^{2\lambda T[g](0)}|x|^{2\nu}(x^2+1)^{(1-\lambda)(\mu-\nu)}e^{2\lambda u(x)}}$$

We rewrite
$$F'(x) = \begin{cases} \dfrac{\kappa e^{\mathfrak{R}_1}}{e^{2\lambda T[g](0)}x^{2\nu}e^{2\lambda u(x)}}, & 0 < x \leq 1 \\ \dfrac{\kappa e^{\mathfrak{R}_2}}{e^{2\lambda T[g](0)}x^{2\nu+2(1-\lambda)(\mu-\nu)}e^{2\lambda u(x)}}, & x \geq 1 \end{cases} \tag{5.14}$$
where $\mathfrak{R}_i$ are bounded above and below by constants depending only on $\mu, \nu$.

From (5.13) we know $0 \leq xu'(x) \leq \mu - \nu$. Since $0 < \lambda \leq \frac{1-2\nu}{2(\mu-\nu)}$, we have $1 - 2\nu - 2\lambda xu'(x) \geq 0$. Therefore
$$\left(\frac{x}{x^{2\nu}e^{2\lambda u(x)}}\right)' = \frac{1}{x^{2\nu}e^{2\lambda u(x)}}(1 - 2\nu - 2\lambda xu'(x)) \geq 0$$

and for $0 < x \leq 1$,
$$\frac{2}{\pi}\frac{\kappa e^{\mathfrak{R}_1}}{e^{2\lambda T[g](0)}}\left(\frac{x}{x^{2\nu}e^{2\lambda u(x)}}\right)' = \frac{2}{\pi}F'(x)(1 - 2\nu - 2\lambda xu'(x)) \tag{5.15}$$
$$= (xu'(x))'(1 - 2\nu - 2\lambda xu'(x))$$

Let $0 < y \leq 1$. Integrating both sides of (5.15) from 0 to $y$ we get
$$\frac{\kappa e^{\mathfrak{R}}}{e^{2\lambda T[g](0)}}\frac{y}{y^{2\nu}e^{2\lambda u(y)}} = yu'(y)(1 - 2\nu - \lambda yu'(y))$$

But
$$\frac{1-2\nu}{2} \leq 1 - 2\nu - \lambda yu'(y) \leq 1 - 2\nu$$

So we have
$$yu'(y) = \frac{\kappa e^{\mathfrak{R}}}{e^{2\lambda T[g](0)}}\frac{y}{y^{2\nu}e^{2\lambda u(y)}} \qquad \text{for } 0 < y \leq 1$$

and
$$e^{2\lambda u(y)}u'(y) = \frac{\kappa e^{\mathfrak{R}}}{e^{2\lambda T[g](0)}}y^{-2\nu}.$$

By integrating from 0 to $x$ we get
$$e^{2\lambda u(x)} - 1 = 2\lambda\frac{\kappa e^{\mathfrak{R}}}{e^{2\lambda T[g](0)}}x^{1-2\nu} \qquad \text{for } 0 < x \leq 1 \tag{5.16}$$

Substitute (5.16) into (5.14), we obtain
$$F'(x) = \frac{\kappa e^{\mathfrak{R}_1}x^{-2\nu}}{e^{2\lambda T[g](0)}(1 + 2\lambda\dfrac{\kappa e^{\mathfrak{R}}}{e^{2\lambda T[g](0)}}x^{1-2\nu})}, \qquad \text{for } 0 < x \leq 1$$

Therefore for $0 < x \leq 1$,
$$F(x) - F(0) = \frac{1}{2\lambda}e^{\mathfrak{R}}\ln(1 + 2\lambda\frac{\kappa e^{\mathfrak{R}}}{e^{2\lambda T[g](0)}}x^{1-2\nu}).$$

This implies for $0 < x \leq 1$

$$0 < \frac{\kappa e^{\Re}}{e^{2\lambda T[g](0)}} x^{1-2\nu} = \frac{1}{2\lambda}(e^{2\lambda(F(x)-F(0))e^{-\Re}} - 1) \leq c(\mu, \nu)$$

for some constant $c(\mu, \nu) > 0$ depending only on $\mu, \nu$. Going back to (5.16), we obtain

$$0 \leq u(x) \leq \frac{1}{2\lambda} \ln(1 + 2\lambda c(\mu, \nu)) \leq c(\mu, \nu) \qquad \text{for } 0 < x \leq 1. \tag{5.17}$$

For $x \geq 1$ the discussion is similar. Notice for $0 < \lambda \leq \frac{2\mu-1}{2(\mu-\nu)}$, $1 - 2\nu - 2(1-\lambda)(\mu-\nu) - 2\lambda x u'(x) \leq 1 - 2\nu - 2(1-\lambda)(\mu-\nu) \leq 0$. Therefore

$$\Big(\frac{x}{x^{2\nu+2(1-\lambda)(\mu-\nu)}e^{2\lambda u(x)}}\Big)'$$
$$= \frac{1}{x^{2\nu+2(1-\lambda)(\mu-\nu)}e^{2\lambda u(x)}}(1 - 2\nu - 2(1-\lambda)(\mu-\nu) - 2\lambda x u'(x)) \leq 0$$

and for $x \geq 1$,

$$\frac{2}{\pi} \frac{\kappa e^{\Re_2}}{e^{2\lambda T[g](0)}} \Big(\frac{x}{x^{2\nu+2(1-\lambda)(\mu-\nu)}e^{2\lambda u(x)}}\Big)'$$
$$= \frac{2}{\pi} F'(x)(1 - 2\nu - 2(1-\lambda)(\mu-\nu) - 2\lambda x u'(x)) \tag{5.18}$$
$$= (xu'(x))'(1 - 2\nu - 2(1-\lambda)(\mu-\nu) - 2\lambda x u'(x))$$

Integrating (5.18) from $\infty$ to $y$, and notice that

$$\lim_{x \to \infty} x u'(x) = \mu - \nu, \qquad \lim_{x \to \infty} \frac{x}{x^{2\nu+2(1-\lambda)(\mu-\nu)}e^{2\lambda u(x)}} = 0,$$

we get for $y \geq 1$,

$$\frac{\kappa e^{\Re}}{e^{2\lambda T[g](0)}} \frac{y}{y^{2\nu+2(1-\lambda)(\mu-\nu)}e^{2\lambda u(y)}}$$
$$= (yu'(y) - \mu + \nu)\{(1 - 2\nu - 2(1-\lambda)(\mu-\nu)) - \lambda(yu'(y) + \mu - \nu)\}$$

Now

$$1 - 2\mu \leq 1 - 2\nu - 2(1-\lambda)(\mu-\nu)) - \lambda(yu'(y) + \mu - \nu)$$
$$= 1 - 2\mu + \lambda(\mu - \nu) - \lambda y u'(y) \leq \frac{1}{2}(1 - 2\mu)$$

Therefore

$$\mu - \nu - yu'(y) = \frac{\kappa e^{\Re}}{e^{2\lambda T[g](0)}} \frac{y}{y^{2\nu+2(1-\lambda)(\mu-\nu)}e^{2\lambda u(y)}}$$

and

$$e^{2\lambda u(y) - 2\lambda(\mu-\nu)\ln y}\Big(\frac{\mu - \nu}{y} - u'(y)\Big) = \frac{\kappa e^{\Re}}{e^{2\lambda T[g](0)}} y^{-2\mu} \tag{5.19}$$

Now $xu'(x) - \mu + \nu \leq 0$, so $u(x) - (\mu - \nu)\ln x$ is decreasing on $[0, \infty)$. Let

$$\lim_{x \to \infty}(u(x) - (\mu - \nu)\ln x) = U_\infty$$

$U_\infty$ is either finite or $-\infty$. If we recall (5.11), Lemma 5.7 and that $\lim_{x \to \infty} T[g](x) = 0$,[11] we know $U_\infty$ must be finite and in fact

$$-2(\mu - \nu) \leq U_\infty + T[g](0) \leq 2(\mu - \nu)$$

---

[11]Because $T[g] \in X$.

Integration (5.19) from $x \geq 1$ to $\infty$, we obtain

$$e^{2\lambda u(x) - 2\lambda(\mu-\nu)\ln x} - e^{2\lambda U_\infty} = 2\lambda \frac{\kappa e^{\Re}}{e^{2\lambda T[g](0)}} x^{1-2\mu} \qquad (5.20)$$

Substitute into (5.14), we get for $x \geq 1$,

$$F'(x) = \frac{\kappa e^{\Re_2} x^{-2\mu}}{e^{2\lambda T[g](0)}\left(e^{2\lambda U_\infty} + 2\lambda \frac{\kappa e^{\Re}}{e^{2\lambda T[g](0)}} x^{1-2\mu}\right)}$$

Integrating from $y \geq 1$ to $\infty$, we obtain

$$F(\infty) - F(y) = \frac{1}{2\lambda} e^{\Re} \ln\left(1 + 2\lambda \frac{\kappa e^{\Re}}{e^{2\lambda(T[g](0)+U_\infty)}} y^{1-2\mu}\right)$$

Therefore for $x \geq 1$,

$$0 \leq \frac{\kappa e^{\Re}}{e^{2\lambda(T[g](0)+U_\infty)}} x^{1-2\mu} = \frac{1}{2\lambda}(e^{2\lambda(F(\infty)-F(y))e^{-\Re}} - 1) \leq c(\mu,\nu)$$

for some constant $c(\mu,\nu)$ depending only on $\mu,\nu$. Going back to (5.20), we have

$$0 \leq u(x) - (\mu-\nu)\ln x - U_\infty \leq \frac{1}{2\lambda}\ln(1 + 2\lambda c(\mu,\nu)) \leq c(\mu,\nu) \qquad \text{for } x \geq 1 \qquad (5.21)$$

Now from (5.17), we know $0 \leq u(1) \leq c(\mu,\nu)$, and from (5.21), $0 \leq u(1) - U_\infty \leq c(\mu,\nu)$. Then $U_\infty$, consequently $T[g](0)$ must be bounded from above and below by two constants that depend only on $\mu,\nu$. Combining (5.17) and (5.21), we have

$$\|u - \frac{1}{2}(\mu-\nu)\ln(x^2+1)\|_{L^\infty[0,\infty)} \leq c(\mu,\nu)$$

for a constant $c(\mu,\nu)$ depending only on $\mu,\nu$. Now from (5.12), because $g \in X_1$ is even, we get

$$\|g\|_\infty \leq c(\mu,\nu)$$

where $c(\mu,\nu)$ is a constant depending only on $\mu,\nu$. This proves Proposition 5.8.

□

We can now conclude from Schaefer's fixed point Theorem that the system (3.12) has a solution $g \in X_1$. This proves the first part of Theorem 3.3. To prove that the solution $g$ satisfies $xg' \in L^\infty(R)$, we notice that from the equation $g = H(F - (\mu-\nu)\arctan\frac{-1}{x})$, where $F$ is as defined in (5.2), we have

$$xg' = H(xF' - (\mu-\nu)\frac{x}{x^2+1}), \qquad (xg')' = H((xF' - (\mu-\nu)\frac{x}{x^2+1})')$$

Using the $L^p$ Boundedness of the Hilbert transform $H$, it is quite straightforward to check that $xF' - (\mu-\nu)\frac{x}{x^2+1} \in L^r(R)$ for $r > \frac{1}{2\mu-1}$, and $(xF' - (\mu-\nu)\frac{x}{x^2+1})' \in L^q(R)$ for $1 < q < \frac{1}{2\nu}$. Therefore we have $xg' \in Y(R)$. This implies $xg' \in L^\infty(R)$, and finishes the proof of Theorem 3.3.

The proof for the a priori estimate in Theorem 3.1 is very much the same as for Proposition 5.8. We just give some main steps. First let $g \in X$, $F$ be defined as in (5.2), and $G(x) = g(x) + \frac{1}{2}(\mu - \nu)\ln(x^2 + 1)$. Assume $g$ is a fixed point of the mapping $T$ defined by (5.1)-(5.2). Then

$$G(x) = G(0) + HF(x) - HF(0)$$

Assume $G$ even, and increasing on $[0, \infty)$. It is easy to check that $F$ satisfies the assumption of Lemma 5.7, therefore in fact $g \in X_1$, and (5.4) holds. Moreover we have from Lemma 5.7,

$$HF(x) - HF(0) = u(x) + \mathfrak{R} \quad \text{for } x > 0$$

where

$$u(x) = \frac{1}{\pi}\int_0^x \frac{F(y) - F(-y)}{y}\,dy$$

and for $x > 0$,

$$xu'(x) = \frac{1}{\pi}(F(x) - F(-x)), \qquad (xu'(x))' = \frac{2}{\pi}F'(x)$$

where

$$F'(x) = \frac{\kappa e^{\mathfrak{R}}}{x^{2\nu} e^{2G(x)}} = \frac{\kappa e^{\mathfrak{R}}}{e^{2G(0)}} \frac{1}{x^{2\nu} e^{2u(x)}}$$

Now

$$\left(\frac{x}{x^{2\nu} e^{2u(x)}}\right)' = \frac{1}{x^{2\nu} e^{2u(x)}}(1 - 2\nu - 2xu'(x))$$

Notice that $xu'(x) = \frac{1}{\pi}(F(x) - F(-x))$ is strictly increasing for $x > 0$, $xu'(x)\big|_{x=0} = 0$, $\lim_{x\to\infty} xu'(x) = \mu - \nu > 1/2 - \nu$. Let $x_0 > 0$ be such that $2x_0 u'(x_0) = 1 - 2\nu$. Then

$$1 - 2\nu - 2xu'(x) > 0, \quad \text{for } 0 < x < x_0, \qquad 1 - 2\nu - 2xu'(x) < 0, \quad \text{for } x > x_0.$$

Now follow very much the same argument as in the proof of Proposition 5.8 on $0 < x \leq x_0$ and $x \geq x_0$. We arrive at the result of Theorem 3.1.

6. Acknowledgement

The author would like to thank Prof. Bo Li and Prof. Lei Ni for their invitation to visit San Diego, CA and for their hospitality. During this visit, the author were able to observe and confirm her analysis with the wave phenomenon. During the past few years, the author have had discussions with friends, colleagues, and many other scientists on waves, singularities and other issues. The current work is influenced and benefited from these conversations. The author started considering self-similar solutions for the water wave equations during her visit to the IMA in 2009-10. The author would like to thank the IMA for its hospitality and financial support.

22  SIJUE WU
[34] P. Zhang, Z. Zhang *On the free boundary problem of 3-D incompressible Euler equations.* Comm. Pure. Appl. Math. V. 61. no.7 (2008), pp. 877-940

University of Michigan, Ann Arbor, MI